\documentclass [11pt] {article}

\usepackage{graphicx}
\usepackage{url}

\begin{document}

\title{Topical Bias in Generalist Mathematics Journals\footnote {Please read and cite the corrected, published article in \textit {Notices of the AMS}, 57(11):1421--1424, 2010.}}

\author{Joseph F. Grcar\thanks{6059 Castlebrook Drive, Castro Valley, CA 94552 USA (jfgrcar@comcast.net).}}

\date{}

\maketitle

Generalist mathematics journals exhibit bias toward the branches of mathematics by publishing articles about some subjects in quantities far disproportionate to the production of papers in those areas within all of mathematics. \textit {Bias\/} is used here because it is the shortest English word with Webster's meaning of ``a tendency of a statistical estimate to deviate in one direction from a true value.'' This paper quantifies the bias, which seems not to be discussed previously, and suggests some consequences of it.

The mathematical topics that are generally agreed to be the major branches of mathematics and the production of papers about them can be determined from the Mathematical Subject Classification and from two databases based on the classification \cite {Fairweather2009}. The classification dates in some form to 1931, when \textit {Zentralblatt f\"ur Mathematik und ihre Grenzgebiete\/} began publishing annual reviews of papers grouped into broad subject areas; \textit {Mathematical Reviews\/} started publishing similar material in 1940. The American Mathematical Society created hierarchical codes to classify papers for its defunct Mathematical Offprint Service in the late 1960s \cite {LeVeque1988}. This AMS (MOS) Classification quickly became the de facto index for mathematical literature. Both \textit {Zentralblatt\/} and \textit {Mathematical Reviews\/} maintain historical databases of indexed papers that can be accessed now through the world wide web. 

The data presented here for the decade 2000--2009 were gathered from the \textit {Zentralblatt\/} database in January, 2010. The 854,547 items that were available at the time will increase as \textit {Zentralblatt\/} completely assimilates papers from the recent past. All the records can be retrieved by Mathematical Subject Classification. 

\begin {figure}[p]
\centering
\includegraphics [scale=0.75] {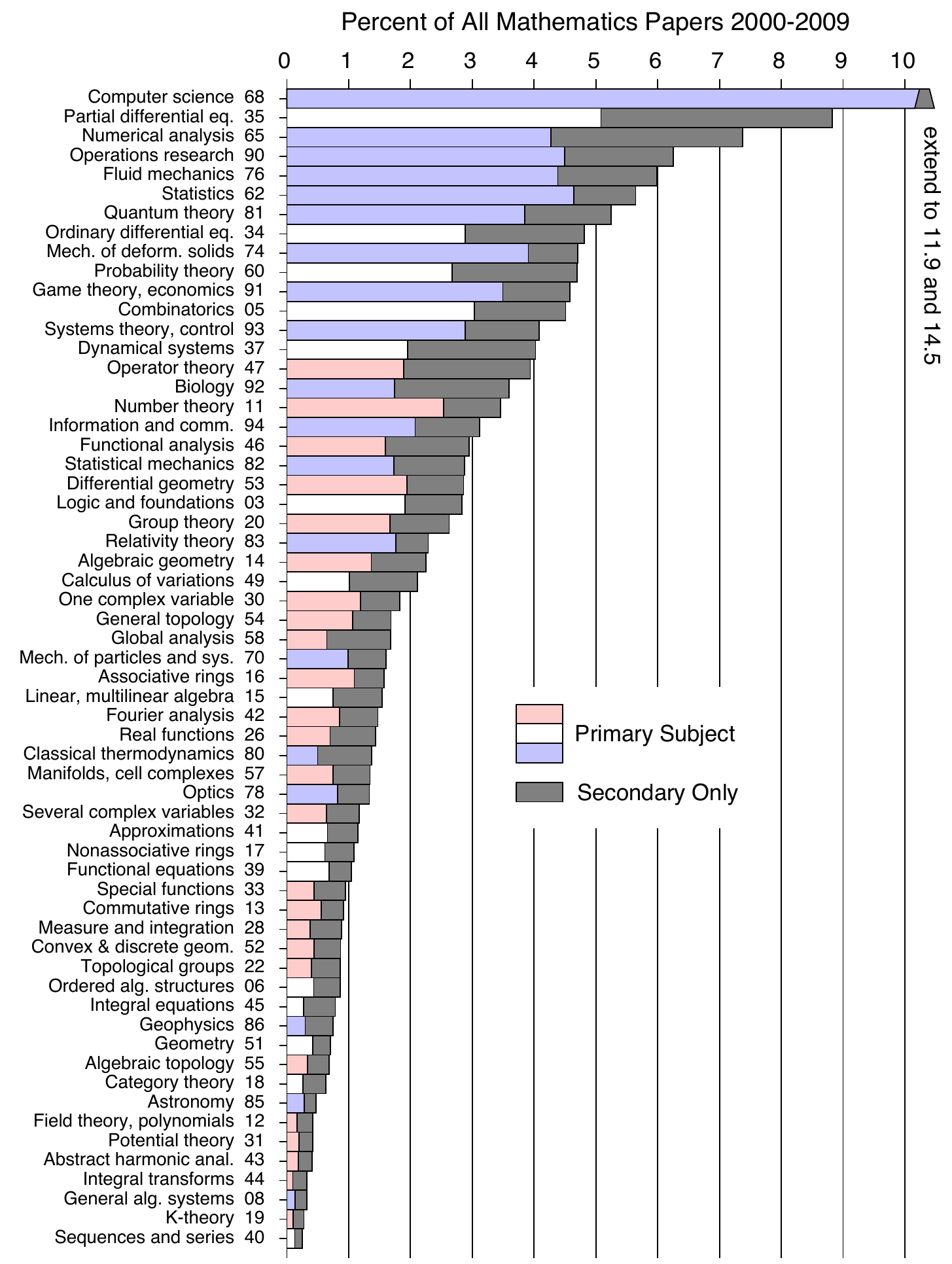}
\caption {The distribution of effort in mathematical research is indicated by the percent of all mathematics papers addressing a particular Mathematical Subject Classification. Three subject classes are omitted: general 00, history 01, and education 97. Papers for which a subject is primary are indicated by red, white, or blue shading (the shading is explained elsewhere). Papers for which a subject is only secondary are indicated by grey shading.}
\label {fig:one}
\end {figure}

Figure \ref {fig:one} displays the percent of all mathematics papers that address each of the major branches of mathematics enumerated by Subject Classification. The present 5-digit codes begin with the 2-digit numbers that reflect the coarsest level of differentiation, namely, the major branches of mathematics. Currently, sixty-three of the 2-digit numbers are assigned. Each paper so catalogued receives one primary code and optionally any number of secondary codes. For example, when the data were gathered, 22,443 papers from 2000--2009 in the \textit {Zentralblatt\/} database had either a primary or a secondary code in the ``group theory'' classification, 20. Thus, approximately $22443 / 854547$ or 2.63 percent of all mathematics papers discussed group theory during 2000--2009. This value is recorded in figure \ref {fig:one}. 

The grey shading in figure \ref {fig:one} indicates papers for which the associated classification is not the primary subject. Such papers are prevalent across mathematics, so they are included in the count of papers for each subject. Since a paper may contribute to several subjects, the percentages for all classes sum to 156.5. The red, white, and blue shades in the figure distinguish subjects in a manner to be explained.

The American Mathematical Society publishes three print journals ``devoted to research articles in all areas of mathematics.'' The extent to which these journals fulfill the Society's promise of inclusiveness can be examined by calculating percentages similar to those in figure \ref {fig:one} but specific to the journals in question. The data presented here are for \textit {Proceedings of the American Mathematical Society}; similar observations can be made for the \textit {Journal\/} and the \textit {Transactions\/}. The \textit {Zentralblatt\/} database holds 4,758 papers from the \textit {Proceedings\/} during the past decade, of which 309 had a primary or secondary code beginning with 20. Therefore, of all papers in the \textit {Proceedings\/}, roughly $309 / 4758$ or $6.49$ percent discuss group theory. This percentage for the generalist journal, $6.49$, contrasts with the $2.63$ percent among all mathematics papers. The fraction of papers about group theory in the journal is over twice the fraction in mathematics as a whole.

The bias of a journal for or against a given branch of mathematics may be defined as the ratio of the fraction of papers about the subject in the journal to the fraction of papers about the subject in all of mathematics. It is convenient to represent the ratio as a base $2$ logarithm, so a bias in favor has a positive value and a bias against has a negative value. The \textit {Proceedings of the American Mathematical Society} thus has a positive bias for group theory of $6.49 / 2.63$ or $2.47 = 2^{+1.31}$. 

\begin {figure}[p]
\centering
\includegraphics [scale=0.75] {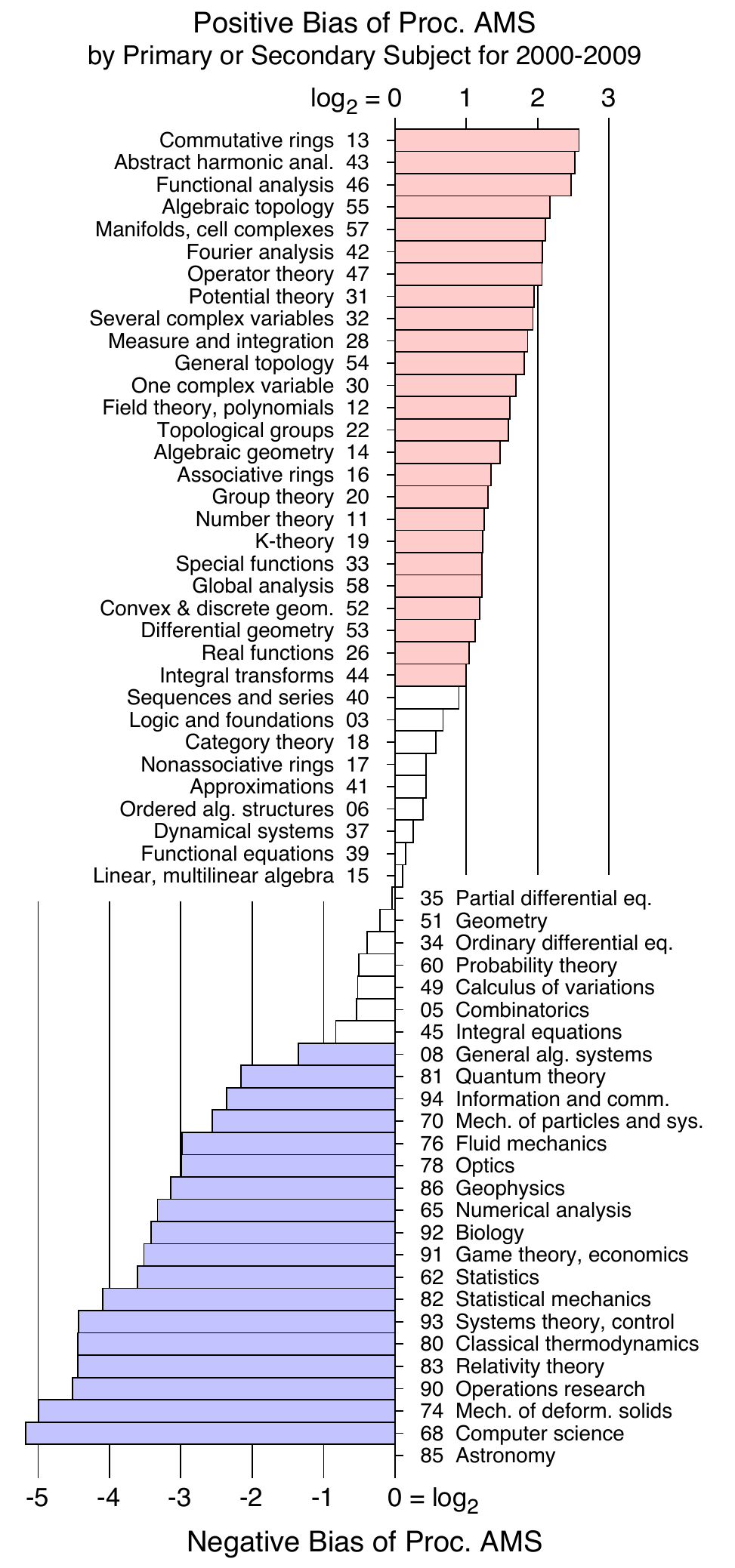}
\caption {Topical bias in {\it Proceedings of the American Mathematical Society\/}. Bias is the ratio of the fraction of publications in the journal for a given subject to the fraction of publications in all of mathematics for that subject. Note the scale is logarithmic to the base 2. The journal had no papers about the one subject at bottom. Subjects with strong positive or negative bias have red or blue shading, respectively, to ease their identification in other figures.}
\label {fig:two}
\end {figure}

Figure \ref {fig:two} displays the biases of the \textit {Proceedings\/} toward all the branches of mathematics. The wide range of values indicates that the journal is unrepresentative of mathematics research. It has a strong bias ($2^{+1}$ to over $2^{+2}$) in favor of twenty-five subjects that are colored red in figures \ref {fig:one} and \ref {fig:two}. Comparing the figures reveals a strong positive bias for three subjects that are of interest to relatively few mathematicians, in that each subject accounts for less than $1$ percent of all mathematics papers:
commutative rings $13$ (bias $2^{+2.58}$),
abstract harmonic analysis $43$ (bias $2^{+2.52}$),
and
algebraic topology $55$ (bias $2^{+2.17}$).
In contrast, the \textit {Proceedings\/} has a neutral to slightly negative bias for three subjects that are of interest to many more mathematicians, in that each of them constitutes over $4$ percent of all mathematics papers:
ordinary differential equations 34 (bias $2^{-0.39}$),
probability theory 60 (bias $2^{-0.51}$),
and
combinatorics 05 (bias $2^{-0.54}$).

The ratio of biases for two subjects is easily seen to equal the ratio of the conditional probabilities that papers about the respective subjects would appear in the journal. Viewed in this way, the biases of the \textit {Proceedings\/} are remarkable. For example, the journal is roughly $6.5 \approx 2^{+2.17} / \, 2^{-0.54}$ times more likely to be the publisher of a paper about algebraic topology than about combinatorics. 

The journal has a strong negative bias ($2^{-1}$ to much below) against eighteen subjects that are colored blue in figure \ref {fig:two}. The same shading in figure \ref {fig:one} reveals the strong negative bias occurs for many branches of mathematics about which most papers are written. Indeed, the \textit {Proceedings\/} is strongly biased against seven of the ten most heavily published subjects. Consequently, the generalist journal neglects subjects to which many mathematicians contribute, or whose very creation is associated with some mathematicians. To cite a few examples:
Tukey \cite {Brillinger2002} for statistics 62, 
Lax \cite {NorwegianAcademy2005} for numerical analysis 65,
Moser \cite {Mather2000} for mechanics 70,
Ladyzhenskaya \cite {Friedlander2005} for fluid mechanics 76,
Dobrushin \cite {Minlos1996} for statistical mechanics 82,
Dantzig \cite {Cottle2007} for operations research 90,
Shannon \cite {Golomb2002} for information theory 94,
and of course,
Wiener \cite {Jerison1995} for systems theory 93,
and von Neumann \cite {Leonard1995} for game theory 91.

\begin {figure}[p]
\centering
\includegraphics [scale=0.75] {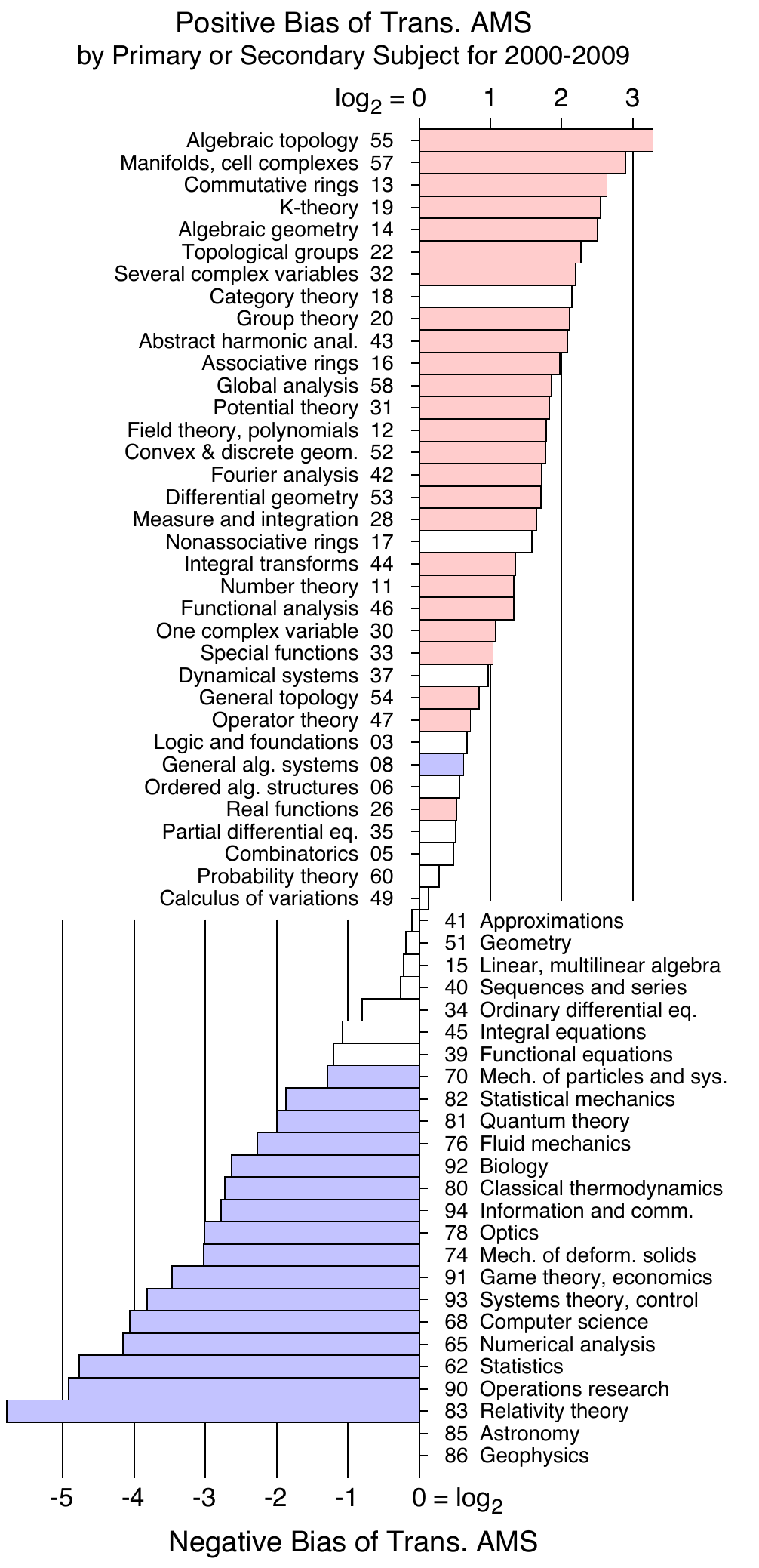}
\caption {Topical bias in {\it Transactions of the American Mathematical Society\/}. The scale is logarithmic to the base $2$. Subjects are shaded as in figure \ref {fig:two} to ease comparison with the {\it Proceedings\/}. The journal had no papers about the two subjects at bottom.}
\label {fig:three}
\end {figure}

The \textit {Transactions\/} has roughly the same biases as the \textit {Proceedings\/} (Figure \ref {fig:three}).  Among heavily published subjects that appear in quantity in either journal, the \textit {Proceedings\/} has stronger positive biases for functional analysis 46 and for operator theory 47.

Explanations for the biases are in a sense circular, in that authors submit where history suggests acceptance is likely, or that editorial boards more accurately evaluate the familiar. For example, not all branches of mathematics employ just the telegraphic style that is prevalent in these journals. Thus biases are self-perpetuating, although surely procedures could be found so submissions that are exceptions to what is usually published do not prove the rule. 

The significant question raised by the data is not how biases occur or how to manage them, but rather, whether the present topical distribution in generalist journals best serves mathematics. Biases in professional journals impart an illusory picture of a field that can be dangerous if it becomes so pervasive as to affect the evolution of the underlying subject matter. In the present case, since the sponsoring Society does not explain the topical selectivity of its generalist journals, and since the literature does not examine the issue, in the absence of clarifying public discussion, readers may assume the contents of the flagship journals confirm prejudices about the branches of the field.  If unchecked over many years, these opinions may influence decisions about curricula, publications, and staffing that can fragment the research community. 

Such dissolution of mathematics may be occurring already as evidenced by the proliferation and growth of fields, particularly since the middle of the last century, that have considerable mathematical content but whose faculties and professional societies have little overlapping membership. When different fields sponsor different branches of mathematics, then the subfields may adopt the cultures of their sponsors, so the branches of mathematics eventually may come to disagree over acceptable idiom, notation, rigor, and terminology. Barriers among the branches of mathematics entail a large opportunity cost, so to speak, because possible collaborations and synergies may not be realized. Moreover, as governments increasingly view scientific research as a component of national wealth, disparities may be expected to grow in the allocation of resources to the fields that encompass different branches of mathematics. In this way the fragmentation of mathematical research cedes to non-mathematicians a greater degree of responsibility to choose which branches of mathematics to encourage and how they should develop. 

\section* {Acknowledgments}  The author is grateful to B.\ Wegner for clarifying the history of the MSC, and to the editor and the referees for advice that improved this article.

\raggedright

\end{document}